\documentclass[10pt]{article}

\usepackage{epsfig}
\usepackage{graphicx}
\usepackage[latin1]{inputenc}
\usepackage{amsmath}
\usepackage{amsthm}
\usepackage{amsfonts}
\usepackage{amssymb}
\usepackage{amstext}
\usepackage{amsgen}
\usepackage{amsbsy}
\usepackage{amsopn}

\usepackage{mathrsfs, a4wide, color}
\usepackage[colorlinks=true, pagebackref]{hyperref}

\numberwithin{equation}{section}

\theoremstyle{plain}
\newtheorem{theorem}{Theorem}[section]

\theoremstyle{definition}

\theoremstyle{remark}

\newcommand\eps{\varepsilon}

\newcommand\vs{\vskip .2cm}

\begin{document}
\title{Wave equations on space-times of low regularity:\\
Existence results and regularity theory in the framework of generalized function algebras}
\author{James~D.E.\ Grant \& Eberhard Mayerhofer\footnote{Faculty of
Mathematics, University of Vienna, Nordbergstrasse 15, 1090
Vienna, Austria, email: \{james.grant, eberhard.mayerhofer\}@univie.ac.at; supported by FWF research grants P16742-N04 and Y237-N13.}}
\date{November 13, 2007}
\maketitle
\thispagestyle{empty}
\begin{abstract}
We present recent developments concerning Lorentzian geometry in algebras of generalized functions. These have, in particular, raised a new interest in refined regularity theory for the wave equation on singular space-times.
\end{abstract}
\medskip
\noindent
\emph{Keywords}: wave equation, generalised Lorentzian geometry, algebras of generalized functions.

\noindent
\emph{MSC 2000}: 35L05, 58J45.

\section{Introduction}
Generalized function algebras have proved themselves a valuable tool
for non-linear analysis of partial differential equations while
maintaining consistency with Schwartz' distributional framework
(\cite{book, Mobook} and the references therein). Their essential
feature is they admit the space of smooth functions as a faithful
subalgebra, which, in view of the Schwartz impossibility
result~\cite{Schw}, is the optimal possibility.

More recent scientific work also focuses on developing a geometric
theory of generalized functions (\cite{CM, Jel} et al.). From the
very beginning of this line of research, one had in mind
applications to fields such as general relativity. The foundations
of generalised pseudo-Riemannian geometry were laid in \cite{gprg},
with further work investigating flows of generalized vector fields
\cite{flows}, and generalised connections on principal
bundles~\cite{connections}.

The work of Vickers and Wilson~\cite{genhyp} on the wave-equation on conical space-times has initiated research on existence and uniqueness results for the initial value problem of the wave-equation for a wider class of singular space-times. A first generalization of their result concerned static-space times \cite{EbMDISS} and, just recently, a quite general existence and uniqueness result has been established, thus dropping the condition of staticity~\cite{mgs}. In the course of this work, it became imperative to study generalized Lorentzian geometry \cite{EbM-AF}, cf.~below.

\vs
The program of this article is to report on recent developments in Lorentzian geometry in a generalised framework and to present results on the wave-equation on singular space-times. Finally, we suggest that regularity information for the wave-equation can be retrieved from refined energy estimates. This study is motivated by concrete examples, cf.~\cite{Jameb}.

\vs
For notation in the field, we refer to the standard reference~\cite{book}.

\section{Lorentzian geometry in algebras of generalized functions}
Loosely speaking, a singular pseudo-Riemannian metric is modelled as a net of smooth metrics,
\[
(g_\eps)_\eps \mod \mathcal N(\mathcal T^0_2(M)),
\]
where $g_\eps$ is a pseudo-Riemannian metric of fixed index for sufficiently small smoothing parameter $\eps$. As generalized functions are not determined pointwise, but only by evaluation on generalized points, one can characterize candidates for generalized pseudo-Riemannian metrics in terms of bilinear forms on the finite dimensional module $\widetilde{\mathbb R^n}$ (for more details, see \cite{gprg}, Proposition 2.\ 1 and \cite{EbM-AF}).
\begin{theorem}\label{chartens02}
Let $g\in \mathcal G^0_2(X)$. The following are equivalent:
\begin{enumerate}
\item \label{chartens021} For each chart $(V_{\alpha},\psi_{\alpha})$ and each $\widetilde x\in (\psi_{\alpha}(V_{\alpha}))^{\sim}_c$ the map
$g_{\alpha}(\widetilde x): \widetilde{\mathbb R}^n\times \widetilde{\mathbb R}^n\rightarrow \widetilde{\mathbb R}$ is symmetric and non-degenerate.
\item \label{chartens022} $g: \mathcal G^0_1(X)\times \mathcal G^0_1(X)\rightarrow \mathcal G(X)$ is symmetric and $\det (g)$ is invertible in $\mathcal G((\psi_{\alpha}(V_{\alpha})))$.
\item \label{chartens023} $\det g$ is invertible in $\mathcal G((\psi_{\alpha}(V_{\alpha})))$ and for each relatively compact open set $V\subset X$ there exists a representative $(g_{\varepsilon})_{\varepsilon}$ of $g$ and $\varepsilon_0>0$ such that $g_{\varepsilon}\mid_V$ is a smooth pseudo-Riemannian metric for all $\varepsilon<\varepsilon_0$.
\end{enumerate}
\end{theorem}
A generalized pseudo Riemannian metric then is defined as a non-degenerate symmetric $g\in \mathcal G^0_2(X)$ with locally constant index. We were motivated by item (\ref{chartens022}) to continue the research in the direction of elaborating a concept of causality from the viewpoint of bilinear forms on $\widetilde{\mathbb R^n}$ on the basis of positivity concepts revisited in the introduction. The first step is, based on a basic perturbation theory, to introduce a well defined index associated to a symmetric matrix $A\in M(n\times n, \widetilde{\mathbb R})$: First, \lq\lq ordered eigenvalues\rq\rq\ of $A$ are defined, $\lambda_1\leq \lambda_2\leq\dots\leq\lambda_n$. Let $j\geq 0$ such that $\lambda_1<<0,\dots,\lambda_j<<0$, $\lambda_{j+1}>>0,\dots,\lambda_{n}>>0$, then $j$ is called the \emph{index\/} of $A$. (This index need not exist.) Finally, the index of a symmetric bilinear form $b$ is then defined as the index of $A:=b(e_i,e_j)$, $e_{i}$ determining the canonical basis of $\widetilde{\mathbb R^n}$. A symmetric bilinear form $b$ is called
\begin{itemize}
\item positive definite, if its index is $0$.
\item Lorentzian, if its index is $1$.
\end{itemize}
Causality is introduced in the following way: Let $b$ be Lorentzian. A vector $v\in\widetilde{\mathbb R^n}$ is called
\begin{itemize}
\item time-like, if $b(v,v)<<0$,
\item null, if $b(v,v)=0$,
\item space-like, if $b(v,v)>>0$
\end{itemize}
Unlike the standard theory, this concept is not a trichotomy (i.e., there are vectors which are neither time-like, null or space-like), since $<<$ does not give rise to a total ordering $<<=$. Essential in the following natural basic considerations is the characterization of free vectors:
\begin{theorem}\label{freechar}
Let $v$ be an element of $\widetilde{\mathbb R}^n$. The following are equivalent:
\begin{enumerate}
\item\label{freechar1} For any positive definite symmetric bilinear form $h$ on $\widetilde{\mathbb R}^n$ we have
\[
h(v,v)>>0
\]
\item\label{freechar2} The coefficients of $v$ with respect to some (hence any) basis span $\widetilde{\mathbb R}$.
\item \label{freechar3}$v$ is free.
\item \label{freeloader6} For each representative $(v_{\varepsilon})_{\varepsilon}\in\mathcal E_M(\mathbb R^n)$ of $v$ there exists some $\varepsilon_0\in I$
such that for each $\varepsilon<\varepsilon_0$ we have $v_{\varepsilon}\neq 0$ in $\mathbb R^n$.
\end{enumerate}
\end{theorem}
Further consequences concern, for instance, characterizations of symmetric positive forms, also, direct summands, the algebraic structure of $\widetilde{\mathbb R}^n$ (being not semi-simple), etc. In order to establish so-called dominant energy conditions for a class of super-energy tensors, it was necessary to establish the inverse Cauchy-Schwarz inequality for time-like vectors in a Lorentzian space-time. We refer the interested reader to the foundational paper \cite{EbM-AF}.

\section{The wave-equation}
The results of the preceding section have been elaborated in close relation with the project of generalising the paper of Vickers and Wilson~\cite{genhyp} for the wave equation on conical space-times. Only the new framework allows to formulate our recent result, as proved in (\cite{mgs}). We give here a brief overview. The initial value problem under consideration corresponds to the scalar wave equation,
\begin{eqnarray}
\Box_{{\bf g}} u &=& f
\nonumber
\\
\label{weqofsettingeps}
u_{\Sigma} &=& v
\\
\nonumber
\widehat{\boldsymbol{\xi}} u_{\Sigma} &=& w,
\end{eqnarray}
where $\Box_{\bf g}$ is the d'Alembertian induced by a generalized Lorentzian metric ${\bf g}$ on a smooth manifold $M$.
Here ${\boldsymbol{\xi}}$ is a $C^\infty$ vector field, and $\widehat{\boldsymbol{\xi}}$ the corresponding unit field; further, $\Sigma$ is the space-like (in the generalized sense) initial surface. The main theorem in our joint paper with R. Steinbauer~\cite{mgs} establishes local existence of scalar solutions to~\eqref{weqofsettingeps}, under the following regularity assumptions (cf.\ also the conclusion section of that paper)
\begin{enumerate}
\item
For all $K$ compact in $M$, for all orders of derivative $k \in \mathbb N_0$ and all $k$-tuples of vector fields $\boldsymbol{\eta}_1, \dots, \boldsymbol{\eta}_k \in \mathfrak X(M)$ and for any representative $( \mathbf{g}_{\eps} )_\eps$ we have:
\begin{itemize}
\item
$\sup_{p \in K}\| \mathscr{L}_{\boldsymbol{\eta}_1} \dots \mathscr{L}_{\boldsymbol{\eta}_k}
\mathbf{g}_{\eps}\|_m=O(\eps^{-k})\quad\quad
(\eps\rightarrow 0)$;
\item
$\sup_{p \in K}\| \mathscr{L}_{\boldsymbol{\eta}_1}\dots\mathscr{L}_{\boldsymbol{\eta}_k}
\mathbf{g}_{\eps}^{-1} \|_m=O(\eps^{-k})\quad\quad
(\eps\rightarrow 0)$.
\end{itemize}

\item
For all $K$ compact in $M$, we have
\[
\sup_{p \in K} \left\Vert \nabla^{\eps} \boldsymbol{\xi}_{\eps} \right\Vert_{\mathbf{m}} = O(\log(\eps)), \quad\quad (\eps \rightarrow 0),
\]
where $\nabla^{\eps}$ denotes the covariant derivative with respect to the Lorentzian metric $\mathbf{g}_{\eps}$.
\item
a causality assumption, to guarantee smooth solutions.
\end{enumerate}

Different arguments suggest that the regularity estimates derived in \cite{mgs} might be improved substantially by refined energy estimates based on the adaption of suitable Sobolev embedding theorems (see, e.g.,~\cite{Aubin, Hebey}) to our setting. In order to apply these results, we require better control on the asymptotic growth of Sobolev-constants corresponding to the regularizing sequences of Riemannian metrics induced on $t = \mathrm{constant}$ slices of our Lorentzian space-time. This, in turn, requires control on the injectivity radius of the induced Riemannian metrics, which is possible if we have control on the diameter, curvature and volume~\cite{CGT}. Investigations on this topic are underway.

\subsection{Regularity issues 1: Wave equations for Lorentzian
metrics of H\"{o}lder regularity}

The class of metrics in the preceding section can be considered as
generalized H\"{o}lder-Zygmund regular of order $s=0$
(cf.~\cite{guezyg} and, for an adaption to the geometric setting, cf.~\cite{Jameb}).
In the latter reference, we have investigated a two
dimensional space-time $(\mathbb{R}^2, \mathbf{g})$, with $\mathbf{g}$
a Lorentzian metric of classic H\"{o}lder regularity $C^{0, \alpha}$,
$\forall\alpha > 0$. The latter is designed in such a way that the
principal symbol of its induced d'Alembertian $\Box_{\mathbf{g}}$ is
precisely the strictly hyperbolic operator used by Colombini \&
Spagnolo~\cite{Colombini} (who present an example of a wave equation
with H\"{o}lder continuous coefficients and smooth data, not admitting
distributional solutions).

In~\cite{Jameb} we bring the metric, by means of a coordinate
transformation $\Phi$ of low regularity $C^{1, \alpha}$ ($\alpha < 1$),
into a form that is evidently conformally flat. This suggests
that solutions $\bar u$ for the wave-equation in flat space-time can
then be pulled back via $\Phi$ to prospective solutions $u$ of
$\Box_{\mathbf{g}}u = 0$. This procedure, however, raises questions
concerning the solution concept, since the coordinate transform is not
$C^2$, the minimal regularity required for the standard analysis to hold.
Nevertheless, for regularisations (that is, $C^\infty$ mappings
$(\Phi_\eps)_\eps$) of $\Phi$, the above is rigorous and we can show
that the smooth pulled back solutions $(u_\eps)_\eps$ even give rise
to generalized solutions, that is, elements of the special algebra
$\mathcal{G}$. So we have, by means of a geometric argument, found
solutions of a non-trivial strictly hyperbolic differential equation
with coefficients of regularity beyond Lipschitz.

\subsection{Regularity issues 2: Collapsing Riemannian manifolds}

An alternative direction for research concerning wave equations on manifolds with singular metrics involves situations where the induced Riemannian metric on three-dimensional hypersurfaces of the Lorentzian manifold \lq\lq collapses\rq\rq\ as $\eps \rightarrow 0$. In particular, in the case of \emph{Cheeger-Gromov collapse\/}~\cite{CG1, CG2}, one constructs nets of Riemannian metrics, $\mathbf{h}_{\eps}$, with the property that, as $\eps \rightarrow 0$, the sectional curvature of $\mathbf{h}_{\eps}$ is uniformly bounded, but the injectivity radius of $\mathbf{h}_{\eps}$ converges uniformly to zero. In this case, the volume of the three-manifold converges to zero as $\eps \rightarrow 0$, but the curvature remains bounded. Taking, for example, a four-dimensional Lorentzian metric of the form $\mathbf{g}_{\eps} = - dt^2 + \mathbf{h}_{\eps}$, we would have asymptotic conditions
\[
\left.
\begin{aligned}
d\mathbf{vol}_{\mathbf{g}_{\eps}} &= O(\epsilon),
\\
\mbox{Curvature} &= O(1)
\end{aligned}
\right\} \qquad \mbox{ as $\eps \rightarrow 0$}.
\label{coeffestimates}
\]
These asymptotic conditions are quite different from those investigated in~\cite{mgs, genhyp} and, as such, will lead to (generalised) solutions of the wave equation that obey quite different asymptotic conditions from those that are standard in Colombeau theory. On the other hand, Cheeger-Gromov collapse is of considerable interest from a differential geometric viewpoint, and is a field where ideas from the field of generalised functions may, perhaps, be usefully applied.


\small
\begin{thebibliography}{10}

\bibitem{Aubin}
{\sc {Aubin, T.}}, {\em Nonlinear analysis on manifolds. {M}onge-{A}mp\`ere
  equations}, vol.~252 of Grundlehren der Mathematischen Wissenschaften,
  Springer-Verlag, New York, 1982.

\bibitem{CG1}
{\sc {Cheeger, J., Gromov, M.}}, {\em Collapsing {R}iemannian manifolds while
  keeping their curvature bounded. {I}}, J. Differential Geom., 23 (1986),
  pp.~309--346.

\bibitem{CG2}
\leavevmode\vrule height 2pt depth -1.6pt width 23pt, {\em Collapsing
  {R}iemannian manifolds while keeping their curvature bounded. {II}}, J.
  Differential Geom., 32 (1990), pp.~269--298.

\bibitem{CGT}
{\sc {Cheeger, J., Gromov, M., Taylor, M.}}, {\em Finite propagation speed,
  kernel estimates for functions of the {L}aplace operator, and the geometry of
  complete {R}iemannian manifolds}, J. Differential Geom., 17 (1982),
  pp.~15--53.

\bibitem{CM}
{\sc {Colombeau, J.~F., Meril, A.}}, {\em Generalized functions and
  multiplication of distributions on {${\mathcal C}^\infty$} manifolds},
  J.~Math.~Anal.~Appl., {\bf 186} (1994), pp.~357--364.

\bibitem{Colombini}
{\sc {Colombini, F. and Spagnolo, S.}}, {\em Some examples of hyperbolic
  equations without local solvability}, Ann. Sci. \'Ecole Norm. Sup. (4), 22
  (1989), pp.~109--125.

\bibitem{Jameb}
{\sc {Grant, J.~D.~E., Mayerhofer, E.}}, {\em Examples of wave-equations on
  singular space-times, preprint}.

\bibitem{mgs}
{\sc {Grant, J.~D.~E., Mayerhofer, E., Steinbauer, R.}}, {\em The wave equation
  on singular space-times}, preprint arXiv:0710.2007v1,  (2007).

\bibitem{book}
{\sc {Grosser, M., Kunzinger, M., Oberguggenberger, M., Steinbauer, R.}}, {\em
  Geometric Theory of Generalized Functions}, vol.~537 of Mathematics and its
  Applications 537, Kluwer Academic Publishers, Dordrecht, 2001.

\bibitem{Hebey}
{\sc {Hebey, E.}}, {\em Nonlinear analysis on manifolds: {S}obolev spaces and
  inequalities}, vol.~5 of Courant Lecture Notes in Mathematics, New York
  University Courant Institute of Mathematical Sciences, New York, 1999.

\bibitem{guezyg}
{\sc {H{\"o}rmann, G.}}, {\em H\"older-{Z}ygmund regularity in algebras of
  generalized functions}, Z. Anal. Anwendungen, 23 (2004), pp.~139--165.

\bibitem{Jel}
{\sc {Jel\'\i nek, J.}}, {\em An intrinsic definition of the {C}olombeau
  generalized functions}, Comment.~Math.~Univ.~Carolinae, {\bf 40} (1999),
  pp.~71--95.

\bibitem{flows}
{\sc {Kunzinger, M., Oberguggenberger, M., Steinbauer, R., Vickers, J.}}, {\em
  Generalized flows and singular {ODE}s on differentiable manifolds}, Acta
  Appl. Math., {\bf 80} (2004), pp.~221--241.

\bibitem{gprg}
{\sc {Kunzinger, M., Steinbauer, R.}}, {\em Generalized {pseudo-}{R}iemannian
  geometry}, Trans. Amer. Math. Soc., {\bf 354} (2002), pp.~4179--4199.

\bibitem{connections}
{\sc {Kunzinger, M., Steinbauer, R., Vickers, J.A.}}, {\em Generalised
  connections and curvature}, Math. Proc. Cambridge Philos. Soc., 139 (2005),
  pp.~497--521.

\bibitem{EbMDISS}
{\sc {Mayerhofer, E.}}, {\em The wave equation on singular space-times}, PhD
  thesis, University of Vienna, 2006.

\bibitem{EbM-AF}
\leavevmode\vrule height 2pt depth -1.6pt width 23pt, {\em On {L}orentz
  geometry in algebras of generalized functions}, Proc. Edinb. Math. Soc., to
  appear,  (preprint arXiv:math-ph/0604052).

\bibitem{Mobook}
{\sc {Oberguggenberger, M.}}, {\em Multiplication of Distributions and
  Applications to Partial Differential Equations}, vol.~{\bf 259} of Pitman
  Research Notes in Mathematics, Longman, Harlow, U.K., 1992.

\bibitem{Schw}
{\sc {Schwartz, L.}}, {\em Sur l'impossibilit\'e de la multiplication des
  distributions}, C.~R.~Acad.~Sci.~Paris, {\bf 239} (1954), pp.~847--848.

\bibitem{genhyp}
{\sc {Vickers, J.~A., Wilson, J.~P.}}, {\em Generalized hyperbolicity in
  conical spacetimes}, Class.~Quantum.~Grav., {\bf 17} (2000), pp.~1333--1360.

\end{thebibliography}
\end{document}